\documentclass{amsart}

\title[Riemann Surfaces]{Embedding subsets of tori Properly into $\CC^2$}
\author{Erlend Forn\ae ss Wold}

\address{Erlend Forn\ae ss Wold \\Department of Mathematics
\\ University of Oslo, P.O. Box 1053, Blindern \\
 0316 Oslo, Norway \\ erlendfw@math.uio.no}

\date{July 26, 2006}

\usepackage{amscd}

\usepackage[left=2.5cm,top=2.5cm,right=2.5cm]{geometry}

\newtheorem{theorem}{Theorem}
\newtheorem{lemma}{Lemma}
\newtheorem{proposition}{Proposition}

\theoremstyle{definition}

\theoremstyle{remark}
\newtheorem{remark}{Remark}

\newcommand{\NN}{\mathbb{N}}
\newcommand{\RR}{\mathbb{R}}
\newcommand{\CC}{\mathbb{C}}

\newcommand{\ZZ}{\mathbb{Z}}
\newcommand{\R}{\mathcal{R}}
\newcommand{\TT}{\mathbb{T}}
\newcommand{\x}{\bf{x}\rm}
\newcommand{\y}{\bf{y}\rm}

\def\d{{\delta}}
\def\e{{\epsilon}}

\def\t{{\tau}}

\def\z{{\zeta}}
\def\l{{\lambda}}
\def\O{{\Omega}}

\usepackage{setspace}

\begin{document}

\large

\bibliographystyle{plain}

\begin{abstract}
Let $\TT$  be a torus.  We prove that all subsets of $\TT$  with
finitely many boundary components (none of them being points)
embed properly into $\CC^2$.  We also show that the algebras of
analytic functions on certain countably connected subsets of
closed Riemann surfaces are doubly generated.
\end{abstract}

\maketitle

\section{Introduction, main results and notation}

Our main concern is the problem of embedding bordered Riemann
surfaces properly into $\CC^2$.  A (finite) bordered Riemann
surface is obtained by removing a finite set of closed disjoint
connected components $D_1,...,D_k$  from a compact surface $\R$,
i.e. the bordered surface is $\tilde\R:=\R\setminus\cup_{i=1}^k
D_i$. \

For a positive integer $d\geq 2$  it is known that there is a
lowest possible integer $N_d=[\frac{3d}{2}]+1$  such that all
Stein manifolds of dimension $d$  embed properly into $\CC^{N_d}$
\cite{eg}\cite{fs2}\cite{sc} (for more details, see for instance
the survey \cite{fc2}). It is also known that all open Riemann
surfaces embed properly into $\CC^3$, but it remains an open
question whether the dimension of the target domain in this case
always can be pushed down to 2. \

For (positive) results when the genus of $\R$  is $0$  we refer to
the articles \cite{kn}\cite{al}\cite{la}\cite{gs}\cite{wd2}, and
in the case of genus $\geq 1$ to \cite{cf}\cite{wd3}. \

We prove the following theorem:

\begin{theorem}\label{main}
Let $\mathbb{T}$  be a torus, and let $U\subset\TT$  be a domain
such that $\TT\setminus U$  consists of a finite number of
connected components, none of them being points.  Then $U$  embeds
properly into $\CC^2$.
\end{theorem}

In \cite{wd3}  we proved that under the assumption that $U$  can
be embedded onto a Runge surface in $\CC^2$, one can embed
arbitrarily small perturbations of $U$  properly into $\CC^2$. Our
task then is to

\

(i) Embed $U$ onto a Runge surface,
\

(ii) Pass from small perturbations to $U$  itself.

\

(We say that a surface $U$  is Runge if holomorphic functions on
$U$  may be approximated uniformly on compacts in $U$  by
polynomials).\

To achieve (i) we recall from \cite{wd3} that for any one
complementary component $D_1$, we have that $\TT\setminus D_1$
embeds into $\CC^2$  by some map $\phi$, and that the image is
Runge.  To embed the smaller domain $U$  onto a Runge surface, we
will perturb the image of $U$ by constructing a map that could be
described as a local (near some neighborhood of $\phi(U)$)
singular shear acting transversally to $\phi(U)$ - the
singularities being placed inside each component of
$\phi(\TT\setminus U)$. This construction is the content of
Section 3. \

To achieve (ii) we will apply a technique from \cite{gs} used by
Globevnik and Stens\o nes to embed planar domains into $\CC^2$. He
and Schramm have shown that any subset of $\TT$ is biholomorphic
to a circular subset $U'$  of another torus $\TT'$ \cite{hs}. This
allows us to identify $U$ with a point in $\RR^N$  in such a way
that the point corresponds to the complex structure on $\TT$ and
the centers and the radii of the boundary components of $U$. Now
small perturbations of $U'$ embeds properly into $\CC^2$, and the
perturbation corresponds to some circled subset of some torus,
i.e. some (other) point in $\RR^N$. So if we identify all subsets
of tori close to $U$  with points in a ball $B$ in $\RR^N$,  we
may in this manner construct a map $\psi:B\rightarrow\RR^N$, such
that all circled domains corresponding to points in the image
$\psi(B)$  embed properly into $\CC^2$.  Our goal is to construct
the map $\psi$  in such a way that it is continuous and close to
the identity.  In that case, by Brouwer's fixed point theorem, the
point corresponding to $U$ will be contained in the image
$\psi(B)$, and the result follows. \

Continuity in the setting of uniformization of subsets of tori is
treated in Section 2, while continuity regarding the
identification of circled subsets with properly embeddable subsets
is dealt with in Section 4. \

As was pointed out in \cite{cf}, the question about the
embeddability of an open Riemann surface $\O$  is related to a
question about the function algebra $\mathcal{O}(\O)$ of all
analytic functions on $\O$.  For an integer $m\in\NN$  we say that
the algebra $\mathcal{O}(\O)$  is $m$-generated if there exist
functions $f_i\in\mathcal{O}(\O)$, $i=1,...,m$  such that
$\CC[f_1,...,f_m]$  is dense in $\mathcal{O}(\O)$.  Since any $\O$
embeds properly into $\CC^3$  we have that $\mathcal{O}(\O)$  is
3-generated, but it is unknown whether or not 2 generators might
be sufficient.  By the perturbation results in Section 3 we get
the following:
\begin{theorem}\label{alg}
Let $\TT$  be a torus, and let $U\subset\TT$ be domain such that
each connected component of $\TT\setminus U$  has got non-empty
interior.  Then the function algebra $\mathcal{O}(U)$ is
$2$-generated.
\end{theorem}
Theorem \ref{alg}  is a special case of the following theorem:
\begin{theorem}\label{algmain}
Let $\R$  be a closed Riemann surface, let $U\subset\R$  be a
domain such that $\partial U$  is a collection of smooth Jordan
curves, and let $\phi:U\rightarrow\CC^2$ be an embedding that
extends across $\partial U$.  Assume that $\phi(\overline U)$ is
polynomially convex. If $V\subset U$ is a connected open set
obtained from $U$  by removing at most countably many disks, then
$\mathcal{O}(V)$ is $2$-generated.
\end{theorem}
The proof of the last two theorems will be given in Section 3. \

As usual we will denote an $\e$-ball centered at a point $p$  in
$\RR^n$ or $\CC^n$  by $B_\e(p)$  (or simply $B_\e$  if the center
is the origin), and the corresponding $\e$-disk in $\CC$  will be
denoted $\triangle_\e(p)$.  By a disk in a Riemann surface $\R$ we
will mean a subset homeomorphic to $\overline\triangle$. \

\verb"Acknowledgement" The author would like to thank Franc
Forstneri\v{c} for several comments and suggestions for
improvements of the present article.  In particular the present
proof of Proposition \ref{perturb} was showed us by
Forstneri\v{c}.

\section{Circled subsets of tori and uniformization}

Let $\t\in\CC$  be contained in the upper half plane $H^+$.  If we
define the lattice
$$
L_\t:=\{m\cdot\t + n\in\CC;m,n\in\ZZ\},
$$
we obtain a torus by considering the quotient $\CC/\sim_\t$, where
$z\sim_\t w\Leftrightarrow z-w\in L_\t$.  It is known that all
tori are obtained in this way.  For a given $\t$  we let
$\R(\Omega(\t))$  denote the quotient, i.e. the torus, and we let
$\Omega(\t)$  denote $\CC$  regarded as its universal cover. We
may choose $\t$  with $0<\mathrm{Re}(\t)\leq 1$. \

We are concerned with subsets of tori with finitely many boundary
components.  Let $\TT$  be a torus, let $\tilde K_1,...,\tilde
K_m$ be compact connected disjoint subsets of $\TT$, such that
$\tilde\TT:=\TT\setminus(\cup_{i=1}^m\tilde K_i)$  is connected.
Then $\TT$ may be identified with its cover $\Omega(\t)$  for some
$\t$, and $\tilde\TT$  with some subset $U$  of $\Omega(\t)$.  It
is clear that $U$  is completely determined by $\t$  and some
choice of complementary components $K_1,\cdot\cdot\cdot, K_m$ of
$U$ that intersect the parallelogram with vertices $0,1,\t,\t+1$,
 and we let $\Omega(\t,K_1,\cdot\cdot\cdot,K_m)$ denote
such a $U$.  We call such a set an m-domain.  We let
$\R(\O(\t,\cdot\cdot\cdot))$  denote the corresponding subset of
$\R(\O(\t))$.  \

Fix an m-domain $\O(\l,K_1,...,K_m)$, and assume that $\l\notin
K_i$  for $i=1,...,m$.  We want to consider a space of m-domains
"close"  to $\O(\l,K_1,...,K_m)$.  For this purpose we recall the
definition of the Hausdorff metric: Let $X$  be a metric space
with distance function $m:X\times X\rightarrow\RR^+$.  For two
closed subsets $S_1,S_2$  of $X$  one defines first
$$
d(S_1,S_2)=\mathrm{sup}_{x\in S_1}\mathrm{inf}\{m(x,y);y\in S_2\}.
$$
Then the Hausdorff distance between the sets $S_1$  and $S_2$  is
defined by
$$
d_H(S_1,S_2)=d(S_1,S_2)+d(S_2,S_1).
$$

Let $\d>0$, let $U_0$  denote the $\d$-disk centered at $\l$, and
for $i=1,...,m$  let $U_i$ denote the $\d$-disk centered at the
closed connected sets $K_i$ with respect to the Hausdorff metric:
$$
U_i=\{S\subset\CC;S \mathrm{\ is \ closed}, d_H(S,K_i)<\d\}.
$$
If $\d$  is small enough then if $\l'\in U_0$  and if $C_i$
 is a connected set $C_i\in U_i$  with $\CC\setminus C_i$ connected for
$i=1,...,m$, then the set $\O(\l',C_1,...,C_m)$ is an m-domain.
(We will also choose $\d$ small enough such that $C_i\in U_i,
C_j\in U_j, i\neq j\Rightarrow C_i\cap C_j=\emptyset$, and such
that no element $C_i\in U_i$ can intersect the disk $U_0$). We
call the set of these m-domains $X^m_{\d}(\O(\l,K_1,...,K_m))$.
Let $\Omega_1=\Omega(\t,K_1,\cdot\cdot\cdot,K_m),\Omega_2=
\Omega(\l,C_1,\cdot\cdot\cdot,C_m)\in
X^m_{\d}(\O(\l,K_1,...,K_m))$, and let $S_1=\{\t\}\cup
K_1\cup\cdot\cdot\cdot\cup K_m$, $S_2=\{\l\}\cup
C_1\cup\cdot\cdot\cdot\cup C_m$  be the corresponding subsets of
$\CC$. We then define
$$
d_1(\Omega_1,\Omega_2):=d_H(S_1,S_2),
$$
As a subset of the set of all m-domains we have all m-domains
whose boundary components are all circles.  We will let these
m-domains be denoted $\Omega(\t,z_1,r_1,\cdot\cdot\cdot,z_m,r_m)$,
where $(z_i,r_i)$ corresponds to the center and the radius of the
ith boundary component (for some choice of ordering of these
components). We will use boldface letters, such as $\bf x\rm$, to
denote a 2m-tuple $\bf x\rm=(z_1,r_1,\cdot\cdot\cdot,z_m,r_m)$ to
simplify notation to $\Omega(\t,\bf x\rm)$.  We call such domains
circled m-domains, and we denote the set of all such domains
$T^m$. \

Let $\Omega(\t,\x)$  be a circled m-domain, and let
$X^m_{\d}(\O(\t,\x))$ be a space as defined above. For all circled
m-domains contained in $X^m_{\d}(\O(\t,\x))$ we have a natural
ordering of all the boundary components, and we may identify all
such domains $\Omega(\l,\y)$ with points $(\l,\y)\in\RR^{2+3m}$.
So if $\e$ is small enough, the points in the ball
$B_\e(\t,\x)\subset\RR^{2+3m}$ are in unique correspondence with
circled m-domains in $X^m_{\d}(\O(\t,\x))$. We may thus give
another metric to this (local) space of circled m-domains,
henceforth denoted $T^m_\e(\t,\x)$, by defining
$$
d_2(\Omega(\t,\x),\Omega(\l,\y)):=\|(\t,\x)-(\l,\y)\|,
$$
where $\|\cdot\|$  is the euclidian distance on $\RR^{2+3m}$. \

We will now give a lemma regarding conformal mappings of arbitrary
m-domains domains onto circular m-domains.  The contents of the
lemma are in essence results proved by He and Schramm \cite{hs}.
Stating the results for the special case of tori, they showed the
following: Let $\TT\setminus\cup_{i=1}^m K_i$ be an m-connected
subdomain of some torus $\TT$. Then there exists some torus $\TT'$
and a domain $\O\subset\TT'$ such that the following holds:

\

(1) $\O$ is circled, meaning that if we lift $\O$  to the
universal cover of $\TT'$  then the complement consists of exact
disks (these disks may also be points), \

(2) $\O$  is conformally equivalent to $\TT\setminus\cup_{i=1}^m
K_i$.

\

Furthermore they proved that

\

(3) A circled domain in the Riemann sphere is unique up to
M\"{o}bius transformations, i.e. if $f:\O_1\rightarrow\O_2$  is a
biholomorphic map between circled domains, then $f$  is the
restriction to $\O_1$  of a M\"{o}bius transformation.

\

Formulating (1) and (2) for m-domains as defined above we have the
following:

\

$(a)$ \ For any $\O=\O(\l,K_1,\cdot\cdot\cdot,K_m)$  there exists
a conformal mapping $f$  that maps $\O$  onto some $\O(\l',\x)\in
T^m$, \

$(b)$ \ The map $f$  respects the relation $\sim_\l$, meaning that
$f(z+m+n\l)=f(z)+m+nf(\l)$  for all $m,n\in\ZZ$.

\

In (b) we have normalized so that $f$ fixes the points $0$ and
$1$.  By $(3)$ we have then that $f$ is unique.

Now fix a domain $\O(\l,K_1,\cdot\cdot\cdot,K_m)$, and consider a
space $X^m_{\d}(\O(\l,K_1,\cdot\cdot\cdot,K_m))$  of nearby
m-domains as defined above.  For each domain
$\O'=\O(\l',C_1,...,C_m)\in
X^m_{\d}(\O(\l,K_1,\cdot\cdot\cdot,K_m))$  there is a unique map
$f$ that maps $\O'$  onto a circular m-domain as above, fixing the
points $0$  and $1$, and we may define a map
$\varphi:X^m_{\d}(\O(\l,K_1,\cdot\cdot\cdot,K_m))\rightarrow T^m$
by
$$
\varphi(\O')=(f(\l'),z_1,r_1,...,z_m,r_m),
$$
where $z_i$  and $r_1$  are the center and radius of the boundary
component corresponding to $C_i$.  Note that by uniqueness, if
$\O'=\O(\l',C_1,...,C_m)$  is a circled m-domain so that $\O'$ has
the representation $\O(\l',z_1,r_1,...,z_m,r_m)$  where
$(z_i,r_i)$ is the center and the radius of $C_i$, then
$\varphi(\O')=(\l',z_1,r_1,...,z_m,r_m)$.  In this respect we may
say that $\varphi\mid_{T^m\cap
X^m_{\d}(\O(\l,K_1,\cdot\cdot\cdot,K_m))} = \mathrm{id}$. \

We will sum these things up in a lemma, and we want to establish
that the map $\varphi$ is continuous.  To prove this we will need
the following definitions and theorem from \cite{gz}:
\

Let $\{B_n\}$, for $n=1,2,..$, denote a sequence of domains in the
Riemann sphere that include the point $z=\infty$.  We define the
kernel of this sequence as the largest domain $B$  including
$z=\infty$ every closed subset of which is contained in each $B_n$
from some $n$ on. We shall say that the sequence $\{B_n\}$
converges to its kernel $B$ if an arbitrary subsequence has the
same kernel $B$. \

\begin{theorem}\label{gz}(\cite{gz}, page 228.)
Let $\{A_n\}$  denote a sequence of domains $A_n$, $n=1,2,...$, in
the Riemann sphere that include the point $z=\infty$.  Suppose
that this sequence converges to a kernel $A$.  Let $\{f_n(z)\}$
denote a sequence of functions $\zeta=f_n(z)$  such that for each
$n=1,2,...$, the function $f_n(z)$  maps the domain $A_n$  onto a
domain $B_n$  including the point $\zeta=\infty$  in such a way
that $f_n(\infty)=\infty$  and $f_n'(\infty)=1$.  Then for the
sequence $\{f_n(z)\}$  to converge uniformly in the interior of
the domain A to a univalent function $f(z)$  it is necessary and
sufficient that the sequence $\{B_n\}$  have a kernel and converge
to it, in which case the function $\zeta=f(z)$  maps A univalently
onto B.
\end{theorem}
We want to apply this theorem for sequences of m-domains.  Let
$A_n$  be a sequence of m-domains including the origin and
converging to an m-domain $A$.  Let $A_n'$  and $A'$  denote the
domains in $\CC$  including $\infty$  given by the correspondence
$z\mapsto\frac{1}{z}$. Then $A_n'$  is a sequence as above, and
$A'$  is its kernel. Let $\{f_n\}$  be a sequence of univalent
functions mapping $A_n$ onto a domain $B_n$  including the origin
and $f_n(0)=0$, $f_n'(0)=1$. For each n define the function
$F_n(z)=\frac{1}{f_n(\frac{1}{z})}$ mapping the domain $A_n'$ onto
$B_n'$, where $B_n'$'s relation with $B_n$  is given by the
correspondence $z\mapsto\frac{1}{z}$. Then the sequences $A_n'$
and $F_n$  satisfy the conditions in the above theorem.  If the
sequence $f_n(z)$  converges to a univalent function $f$  on $A$,
the sequence $F_n$  converges to a univalent function $F$  on
$A'$.  By the theorem the sequence $B_n'$ has a kernel $B'$ and
converges to it, and $F$  maps $A'$  onto $B'$.  This implies that
the sequence $B_n$ has a kernel $B$  and converges to it, and $f$
maps $A$  onto $B$.  On the other hand, if the sequence $B_n$  has
a kernel $B$ and converges to it, then the sequence $B_n'$ has a
kernel $B'$ and converges to it, and by the theorem $F_n$
converges to a univalent function $F$ on $A'$, mapping $A'$ onto
the kernel $B'$. So the sequence $f_n$ converges to a univalent
function $f$ on $A$ mapping $A$  onto the kernel $B$.

\begin{lemma}\label{mmap}
Let $X^m_{\d}(\O(\t,K_1,\cdot\cdot\cdot,K_m))$  be a space of
m-domains as defined above.  There is a map
$\varphi:X^m_{\d}(\O(\t,K_1,\cdot\cdot\cdot,K_m))\rightarrow T^m$
such that the following holds:

\

(i) $\R(\varphi(\O'))$  is conformally equivalent to $\R(\O')$ for
all $\O'\in X^m_{\d}(\O(\t,K_1,\cdot\cdot\cdot,K_m))$, \

(ii) $\varphi\mid_{T^m\cap
X^m_{\d}(\O(\t,K_1,\cdot\cdot\cdot,K_m))} = \mathrm{id}$,\

(iii) $\varphi$  is continuous with respect to $d_1$  and $d_2$. \
\end{lemma}

\begin{proof}
We have already defined $\varphi$  and established $(i)$  and
$(ii)$.  To prove continuity we first choose a different
normalization of the uniformizing maps. For each map
$f:\O'\rightarrow\CC$ as above, we compose with a linear map and
assume that $f(0)=0, f'(0)=1$. \

Let $\O(\l,Y_1,...,Y_m)\in
X^m_{\d}(\O(\t,K_1,\cdot\cdot\cdot,K_m))$  and let
$f:\O(\l,Y_1,...,Y_m)\rightarrow\CC$ be the corresponding map. Let
$\{\O(\l_j,Y_1^j,..,Y_m^j\}\subset
X^m_{\d}(\O(\t,K_1,\cdot\cdot\cdot,K_m))$  such that
$\O(\l_j,Y_1^j,...,Y_m^j)\rightarrow \O(\l,Y_1,...,Y_m)$ and let
$f_j:\O(\l_j,Y_1^j,...,Y_m^j)\rightarrow\CC$ be the corresponding
maps for those domains.  By abuse of notation we will let $f(Y_i)$
and $f_j(Y^j_i)$  denote complementary components of the images.
Note that the sequence of domains $\O(\l_j,Y_1^j,...,Y_m^j)$  has
the domain $\O(\l,Y_1,...,Y_m)$  as its kernel and converges to
it. We claim that $f_j\rightarrow f$ uniformly on compacts in
$\O(\l,Y_1,...,Y_m)$, and that $f_j(Y_i^j)\rightarrow f(Y_i)$.
This will prove the continuity of the map $\varphi$ defined above.
That we chose a different normalization does not matter since we
will then also have that
$\frac{f_j}{f_j(1)}\rightarrow\frac{f}{f(1)}$.  \

To show that $f_j\rightarrow f$  it suffices to show that every
subsequence of $f_j$  admits a subsequence converging to $f$. By
assumption on the family
$X^m_{\d}(\O(\t,K_1,\cdot\cdot\cdot,K_m))$ there exists a $t>0$
such that $\overline\triangle_t=\{\z\in\CC;|\z|\leq t\}\subset
\O(\l_j,Y_1^j,...,Y_m^j)$ for all $j$.  Now let $t_0<t$  and
consider the functions $h_j(z)=\frac{1}{f_j(z)}$  on
$W^j_{t_0}=\O(\l_j,Y_1^j,...,Y_m^j)\setminus\overline{\triangle}_{t_0}$.
By Koebe's $\frac{1}{4}$-Theorem we have that
$h_j(W^j_{t_0})\subset\triangle_{\frac{4}{t_0}}$  for all $j$, so
the sequence $h_j$  is a normal family on
$W_{t_0}=\O(\l,Y_1,...,Y_m)\setminus\overline{\triangle}_{t_0}$.
Passing to a subsequence we assume that $h_j\rightarrow h$. Now
$h$ cannot be constantly zero, for this would mean that
$f_j=\frac{1}{h_j}\rightarrow\infty$ uniformly on compacts. This
would contradict the fact that $f_j'(0)=1$ for all $j$. But this
means that that the sequence $f_j$ converges to some function $g$
on $W_{t_0}$, hence we may assume that $f_j$ converges to $g$  on
$\O(\l,Y_1,...,Y_m)$. Since $g'(0)=1$ we have that $g$  cannot be
constant, and we conclude that $g$ maps $\O(\l,Y_1,...,Y_m)$
univalently onto some subset of $\CC$. \

Since $f_j$  converges to $g$  we have now that the for each $i$,
the set $f_j(Y^j_i)$  is a bounded sequence of disks
$\triangle_{r^j_i}(z^j_i)$  (some of these disks could be points).
So by passing to a subsequence we may assume that each of the
sequence of pairs $(z^j_i,r^j_i)$ converges to some pair
$(z_i,r_i)$. We have that

\

$(B) \ f_j(z+m+n\l_j)=f_j(z)+m f_j(1)+n f_j(\l_j)$

\

for all $j$ and for all $m,n\in\ZZ$.  So if we let $Q_j$  be the
set of disks in $\CC$  generated by the set of disks
$\triangle_{r^j_i}(z^j_i)$  and the lattice determined by
$f_j(1)$  and $f_j(\l_j)$, we get that
$f_j(\O(\l_j,Y^j_1,...,Y^j_m))=\CC\setminus Q_j$. \

From $(B)$  we now get that

\

$(C) \ g(z+m+n\l)=g(z)+mg(1)+ng(\l)$

\

for all $m,n\in\ZZ$. \

We must have that $g(1)$  and $g(\l)$  are linearly independent
over $\RR$.  To see this let $V$  be some open set in
$\O(\l,K_1,...,K_m)$ containing the point $\l$.  Then $g(V)$
contains an open set around $g(\l)$.  Now for each $m,n\in\ZZ$ let
$V_{m,n}$  denote the translated sets $V+m+n\l$. Then
$g(V_{m,n})=g(V)+mg(1)+ng(\l)$, and if $g(1)$  and $g(\l)$ are
linearly dependent over $\RR$  then $g(V_{m,n})$  would intersect
the straight line segment between $0$  and $g(\l)$  for infinitely
many choices of $m,n\in\ZZ$.  This would contradict the fact that
$g$  is univalent. \

Let $Q$  now denote the circled subset of $\CC$  generated by the
disks $\triangle_{r_i}(z_i)$  and the lattice determined by $g(1)$
and $g(\l)$.  Now $\CC\setminus Q$  is the kernel for sequence
$\CC\setminus Q_j$, and it follows from Theorem \ref{gz}  that
$g(\O(\l,Y_1,...,Y_m))=\CC\setminus Q$.  But then $g$  is the
unique function satisfying $g(0)=0, g'(0)=1$  that maps
$\O(\l,Y_1,...,Y_m)$ onto a circled subset of $\CC$  having a
cluster point at infinity, and this contradicts $(A)$.  We
conclude then that $f_j\rightarrow f$.
\

Now from Theorem \ref{gz}  we have that $f(\O(\l,Y_1,...,Y_m))$ is
the kernel for the sequence $f_j(\O(\l,Y^j_1,...,Y^j_m))$  to
which it converges.  Since an arbitrary subsequence has the same
kernel we have that each sequence of disks $f_j(Y^j_i)$  must
converge to $f(Y_i)$, and this completes the proof.
\end{proof}

Now let $\O(\t,\x)\in T^m$  so that no boundary components
intersects the point $\t$, let $X_\d^m(\O(\t,\x))$ be a space as
defined above, and choose $\e>0$ such that $T^m_\e(\t,\x)\subset
X_\d^m(\O(\t,\x))$. Let $\varphi:X_\d^m(\O(\t,\x))\rightarrow T^m$
be the map from Lemma \ref{mmap}.  We then have the following:

\begin{lemma}\label{lim}
For every $\mu>0$ there exists a $\widehat\d>0$ such that, if
$$
\psi:T^m_\e(\t,\x)\rightarrow X_\d^m(\O(\t,\x))
$$
is a map with $d_1(\psi(\O(\l,\y)),\O(\l,\y))<\widehat\d$ for all
$\O(\l,\y)\in T^m_\e(\t,\x)$, then
$$
d_2(\varphi\circ\psi(\O(\l,\y)),\O(\l,\y))<\mu
$$
for all $\O(\l,\y)\in T^m_\e(\t,\x)$.
\end{lemma}

\begin{proof}
This follows from the facts that $\varphi|_{T^m\cap
X_\d^m(\O(\t,\x))}=\mathrm{id}$, $\varphi$ is continuous, and
$\overline{T^m_\e(\t,\x)}$  is complete.
\end{proof}

Theorem \ref{main}  will follow from the previous lemmas and the
following proposition.  The proof of the proposition will be given
in sections 3 and 4.

\begin{proposition}\label{mainmap}
Let $\O(\t,\x)\in T^m$  such that no complementary component of
$\O(\t,\x)\in T^m$  intersect the point $\t$, and such that no
boundary component is a single point.  Let $X^m_\d(\O(\t,\x))$ be
a space as above. If $\e>0$ is small enough, then for all
$\widehat\d>0$ there exists a map
$\psi:T^m_\e(\O(\t,\x))\rightarrow X^m_\d(\O(\t,\x))$ such that
the following holds:

\

(i) $\psi$  is continuous with respect to $d_1$  and $d_2$,
\

(ii) $d_1(\O(\l,\y),\psi(\O(\l,\y)))<\widehat\d$  for all
$\O(\l,\y)\in T^m_\e(\O(\t,\x))$, \

(iii) All $\R(\psi(\O(\l,\y)))$  embed properly into $\CC^2$.
\end{proposition}

\emph{Proof of Theorem} \ref{main}: \ Lift $U$  to the universal
cover of $\TT$  and write this lifting as an m-domain
$\O(\l,K_1,...,K_m)$.  By Lemma 1, $\O(\l,K_1,...,K_m)$ is
biholomophic to some circled m-domain $\O(\t,\x)\in T^m$(see
(1),(2),(a) and (b) on page 3), so it is enough to proof the
result for $\R(\O(\t,\x))$.  By a linear translation we may assume
that no boundary component of $\O(\t,\x)$  intersect the point
$\t$, and we cannot have that any boundary component of
$\O(\t,\x)$  is a point, since no $K_i$  is a point.  Let $\e>0$
be in accordance with Proposition \ref{mainmap}. There exists a
$\mu>0$ such that if $F:B_\e(\t,\x)\rightarrow\RR^{2+3m}$ is a
continuous map satisfying
$$
(*) \ \|F-id\|_{B_\e(\t,\x)}<\mu,
$$
then
$$
(**)(\t,\x)\in F(B_\e(\t,\x)).
$$
Choose $\widehat\d>0$ depending on
$\mu$  as in Lemma \ref{lim}, choose $\psi$  as in Proposition
\ref{mainmap} depending on $\widehat\d$, and consider the
composition
$$
F=\varphi\circ\psi
$$
(regarded as a map from $B_\e(\t,\x)$  into $\RR^{2+3m}$).  Then
$F$ is a map satisfying $(*)$  so we have $(**)$.  We have that
all circled m-domains corresponding to points in $F(B_\e(\t,\x))$
embed properly into $\CC^2$, so $\R(\O(\t,\x))$ embeds properly
into $\CC^2$. $\square$ \

It is clear that we have proved the following formulation of
Theorem 1, which we formulate for easier reference in applications
to embeddings with interpolation: \

\bf{Theorem 1':}\rm \ Let $\TT$  be a torus, and let $U\subset\TT$
be a domain such that $\TT\setminus U$  consists of a finite
number of connected components, none of them being points.  Then
$U$ embeds onto a surface in $\CC^2$  satisfying the conditions in
Theorem 1 in \cite{wd3}.

\section{Perturbing surfaces in $\CC^2$  and consequences for function algebras.}

Let $\mathcal{R}$  be an open Riemann surface, and let $U$  be an
open subset of $\mathcal{R}$.  We say that $U$  is Runge in
$\mathcal{R}$  if every holomorphic function $f\in\mathcal{O}(U)$
can be approximated uniformly on compacts in $U$  by functions
that are holomorphic on $\mathcal{R}$.  If $\phi(\mathcal{R})$  is
an embedded surface in $\CC^2$  we will say that
$\phi(\mathcal{R})$  is Runge (in $\CC^2$)  if all functions
$f\in\mathcal{O}(\phi(\mathcal{R}))$  can be approximated
uniformly on compacts in $\phi(\mathcal{R})$  by polynomials. \
Now let $M$  be a complex manifold and let $K\subset M$  be a
compact subset of M.  Recall the definition of the holomorphically
convex hull of $K$  with respect to $M$:
$$
\widehat{K}_M=\{x\in M;|f(x)|\leq\|f\|_K, \forall
f\in\mathcal{O}(M)\}.
$$
If $M=\CC^n$  we simplify to $\widehat K=\widehat K_{\CC^n}$, and we call
$\widehat K$  the polynomially convex hull of $K$.  If $K=\widehat K$  we say that $K$  is
polynomially convex. \

For an open Riemann surface $\mathcal{R}$, and a compact set
$K\subset\mathcal{R}$, we have that:

\

(1) $\widehat K_{\mathcal{R}}$ is the union of $K$  and all the
relatively compact components of $\mathcal{R}\setminus K$, \

(2) An open subset $U$  of $\mathcal{R}$  is Runge if and only if
$\widehat K_{\mathcal{R}}\subset U$  for all compact $K\subset U$.

\

These results can be found in \cite{bs}, \cite{ma}. \

We will need the following standard result:
\begin{lemma}\label{runge}
Let $U\subset\CC^k$  be Runge and Stein, and let $X\subset U$  be
an analytic set.  For $M\subset\subset X$  we have that
$$
\widehat M=\widehat M_{\mathcal{O}(U)}=\widehat
M_{\mathcal{O}(X)}.
$$

\end{lemma}
\begin{proposition}\label{perturb}
Let $\R$  be a closed Riemann surface, let $V\subset\R$  be a
domain such that $\partial V$  is a collection of smooth Jordan
curves, and let
$$
\phi:V\rightarrow\CC^2
$$
be an embedding, holomorphic across the boundary.  Assume that
$\phi(\overline V)$ is polynomially convex. Then for any finite
set of distinct points $\{p_i\}_{i=1}^m\subset V$, there exist
arbitrarily small open disks $D_i\subset V$  with $p_i\in D_i$,
and a neighborhood $\O$ of $\phi(\overline{V}\setminus\cup_{i=1}^m
D_i)$, such that for all $\e>0$ there exists an injective
holomorphic map
$$
\xi:\O\rightarrow\CC^2
$$
such that the following holds:

\

(i) $\|\xi-id\|_{\phi(\overline{V}\setminus\cup_{i=1}^m D_i)}<\e$
\

(ii) $\xi\circ\phi(\overline{V}\setminus\cup_{i=1}^m D_i)$  is
polynomially convex.
\end{proposition}
\begin{proof}
Let $V\subset\subset W$  such that $\phi|_W$ is an embedding.
Since $\phi(\overline V)$  is polynomially convex there is a Runge
and Stein neighborhood basis $U_j$ of $\phi(\overline V)$ in
$\CC^2$. We may assume that $W_j:=\phi(W)\cap U_j$  is a closed
submanifold of $U_j$  for all $j\in\NN$, and that $\phi(V)$  is
Runge in $W_j$.  Let $x_i$ denote $\phi(p_i)$
 for $i=1,...,m$, and let $Q=\{x_1,...,x_m\}$. \

Now let $\mathcal{N}$  denote the normal bundle of $W_1$.  Since
$\mathcal{N}$  is a line bundle and $W_1$ is a Riemann surface, we
have that $\mathcal{N}\cong W_1\times\CC$ (see for instance
\cite{fs}, p.229).  For some large enough $j\in\NN$  we have that
$U_j$ embeds into $\mathcal{N}$ with $W_j$ as the zero section,
i.e. there is an injective holomorphic map
$$
(*) \ F:U_j\rightarrow W_j\times\CC
$$
such that $F(x)=(x,0)$  for all $x\in W_j$.  We might as well
assume that this is true for $j=1$ (for a reference to these
claims see \cite{gr} pages 255-258 and Remark \ref{outline}
below).\

Let $f\in\mathcal{O}(\phi(W))$ with $f(x)=0$  for $x\in Q$, and
$f(x)\neq 0$ for $x\notin Q$  (see for instance \cite{fs}).  For
any $\d>0$ we let
$$
\psi_\d: W_1\setminus Q\times\CC\rightarrow W_1\setminus
Q\times\CC
$$
be the biholomorfic map defined by
$\psi_\d(x,\l)=(x,\l+\frac{\d}{f(x)})$.  Then
$\psi_\d(F(W_1\setminus Q))$  is a closed submanifold of
$W_1\times\CC$  for all choices of $\d$, and we get that
$W_1^\d:=F^{-1}(\psi_\d(F(W_j\setminus Q)))$  is a closed
submanifold of $U_1$. \

Let $\O_j$  be a neighborhood basis of
$\phi(\overline{V}\setminus(\cup_{i=1}^m D_i))$  in $\CC^2$.  If
$j$ is large enough and $\d$  is small enough we have that
$$
G_\d:=F^{-1}\circ\psi_\d\circ F:\O_j\rightarrow U_1
$$
is an injective holomorphic map.  Moreover we have that
$G_\d(\phi(\overline{V}\setminus(\cup_{i=1}^m D_i)))$ is
holomorhically convex in $W_1^\d$.  Put $\O:=\O_j$, $\xi:=G_\d$,
and the result follows by Lemma \ref{runge}.
\end{proof}
\begin{remark}\label{outline}
We outline a simple proof of the existence of the map $(*)$ in our
setting: Let $g\in\mathcal{O}(U_1)$ be a defining function for
$W_1$, and let $\bigtriangledown g(x)$ denote the gradient of $g$.
Such a function exists since Cousins second problem has a solution
in this setting.  Define a map
$$
H:W_1\times\CC\rightarrow\CC^2
$$
by $H(x,\l)=x +\l\cdot\bigtriangledown g(x)$.  It is seen that
$H$  is injective near $W_1\times\{0\}$, and we may let $F=H^{-1}$
on $U_j$  if $j$  is big enough.
\end{remark}
\emph{Proof of Theorem \ref{algmain}:} Let $\{K_j\}$  be a
holomorphically convex exhaustion of $V$  such that $U\setminus
K_j$  has finitely many complementary components for each
$j\in\NN$. We will repeatedly use Proposition \ref{perturb}  to
construct an embedding $\phi$  of $V$  into $\CC^2$  such that
each $\phi(K_i)$ is polynomially convex, and this will prove the
theorem.\

Assume that we are in the following situation which we call
$S_i$:\

We have found a domain $V_i\subset\R$ such that $V\subset\subset
V_i$, with $K_i$ holomorphically convex in $V_i$, and an embedding
$\phi_i:V_i\rightarrow\CC^2$ such that the conditions in
Proposition \ref{perturb}  are satisfied for the pair
$(V_i,\phi_i)$.  In particular we have that $\phi_i(K_i)$  is
polynomially convex. \

We will show that we can use Proposition \ref{perturb}  to pass to
situation $S_{i+1}$. \

Let $T_1,...,T_k$  denote the connected components of
$V_i\setminus K_{i+1}$.  If no $T_j$  is relatively compact in
$V_i$  we have that $K_{i+1}$  is holomorphically convex in $V_i$
and we define $V_{i+1}:=V_i, \phi_{i+1}:=\phi_i$.  Assume on the
other hand that $T_{i_1},...,T_{i_s}$ are relatively compact in
$V_i$. By assumption and since $K_{i+1}$ is holomorphically convex
in $V$, we may find points $p_j\in T_{i_j}$ such that $p_j\in
(U\setminus V)^\circ$.  And so there are disks $D_j\subset
V_i\setminus V$ such  that $p_j\in D_j$. Define
$V_{i+1}=V_i\setminus\cup_{j=1}^s D_j$  and Proposition
\ref{perturb}  furnishes the map $\phi_{i+1}$.  We are in
$S_{i+1}$. \

We may now use this procedure to construct an appropriate
embedding of $V$ into $\CC^2$.  Let $V_1$  be a smoothly bounded
domain in $\R$, homeomorphic to $U$  with $U\subset\subset V_1$,
and such that $\phi$ is defined on $V_1$.  Assume that $K_1$  is a
point and define $\phi_1:=\phi$.  Notice that for each step, when
passing from $S_i$ to $S_{i+1}$, we may choose any $\d_i>0$  and
make sure that $\|\phi_{i+1}-\phi_i\|_{K_{i+1}}<\d_i$.  Therefore
we may choose a sequence $\{\phi_i\}$  such that
$$
\phi:=\lim_{i\rightarrow\infty}\phi_i
$$
exists on $V$  and is an embedding.  Moreover, since $\phi_i(K_i)$
is polynomially convex for each $i\in\NN$, and since
$\phi(K_{i+1})$ can be made an arbitrarily small perturbation of
$\phi_i(K_{i+1})$, we may assume that each $\phi(K_i)$  is
polynomially convex. The result follows.  $\hfill\square$
\

\emph{Proof of Theorem \ref{alg}:} Let $T$  be a connected
component of $\TT\setminus V$, and let $p\in T$  be an interior
point.  Then $\TT\setminus\{p\}$  embeds as a closed submanifold
of $\CC^2$  by some map $\phi$.  Let $D$  be a smoothly bounded
disk such that $D\subset\subset T$, and define $U=\TT\setminus D$.
The collection $(U,\phi,V)$  satisfies the conditions in Theorem
\ref{algmain}.$\hfill\square$

\section{Continuous perturbation of families of Riemann surfaces -
proof of Proposition \ref{mainmap}}

Briefly the idea behind the proof of Proposition \ref{mainmap} is
the following: Start with the space $T^m_\e(\O(\t,\x))$  and
consider Theorem \ref{old}  below.  In effect we showed in
\cite{wd3} that for each fixed $\O(\l,\y)\in T^m_\e(\O(\t,\x))$
there exists an arbitrarily small perturbation $U_{(\l,\y)}$  of
$\O(\l,\y)$  such that $U_{(\l,\y)}$  embeds onto a surface in
$\CC^2$  satisfying the conditions in Theorem \ref{old}. I.e.
$U_{(\l,\y)}$  embeds properly into $\CC^2$.  Suppose that we
could make sure that the perturbed $m$-domains vary continuously
with the parameter $(\l,\y)$  (with respect to the metric defined
in Section 2). Then the correspondence $\O(\l,\y)\mapsto
U_{(\l,\y)}$ defines a continuous map
$\psi:T^m_\e(\O(\t,\x))\rightarrow X^m(\t,\x)$, and all the image
domains embed properly into $\CC^2$.  If $\psi$ could be made
arbitrarily close to the identity then Proposition \ref{mainmap}
would follow from Lemma \ref{lim}.  This is indeed what we will
prove. \

The following theorem is approximately the same as Theorem 1 in
\cite{wd3}.  The difference is that Theorem 1 was formulated for
surfaces with smooth boundaries, whereas the following is
formulated for surfaces with piecewise smooth boundaries.  The
difference in the proof however is not significant.
\begin{theorem}\label{old}
Let $M\subset\CC^2$  be a Riemann surface whose boundary
components are piecewise smooth Jordan curves
$\partial_1,...,\partial_m$. Assume that there are points
$p_i\in\partial_i$  such that
$$
\pi_1^{-1}(\pi_1(p_i))\cap\overline M=p_i.
$$
Assume that each boundary component $\partial_i$ is smooth near
$p_i$, and that all points $p_{i}$  are regular points of the
projection $\pi_{1}$. Then $M$ can be properly holomorphically
embedded into $\CC^2$.
\end{theorem}
As outlined above we want to embed families of $m$-domains onto
surfaces satisfying the conditions in this theorem.  It seems
worth it however to formulate a result for closed Riemann surfaces
in general: Fix an integer $g\geq 0$.  Let $B_\e$  denote a ball
of radius $r=\e$  in some $\RR^N$  and let $X$  be a smooth
manifold with a projection $\pi:X\rightarrow B_\e$  such that
$X_y:=\pi^{-1}(y)$ is a closed Riemann surface of genus $g$  for
each $y\in B_\e$ $-$ the complex structure on each fibre $Y_y$
being specified by the parameter $y$.  Let $m:X\times
X\rightarrow\RR^+$ be a smooth metric on $X$ that induces the
topology. \

For $i=1,...,m$  let $f_i:B_\e\times\overline\triangle$  be a
smooth embedding  such that
$f_i(\{y\}\times\overline\triangle)\subset X_y$  for each $y\in
B_\e$, and such that the images
$f_i(B_\e\times\overline\triangle)$  are pairwise disjoint.  Let
$Y:=X\setminus\cup_{j=1}^m f_j(B_\e\times\overline\triangle)$.
Then $Y$  is a submanifold of $X$  and each fiber $Y_y\subset X_y$
is an open Riemann surface (specifically a closed Riemann surface
of genus $g$ with $m$  disks removed).  For $0<\d<1$  let $Y^\d$
denote $X\setminus\cup_{j=1}^m
f_i(B_\e\times\overline\triangle_{1-\d})$.
\begin{proposition}\label{conper}
Let $F:Y^\d\rightarrow B_\e\times\CC^2$  be a smooth map such that
$F(y,\cdot):Y^\d_y\rightarrow\{y\}\times\CC^2$  is a holomorphic
embedding for each $y\in B_\e$.  Assume that $F(y,\overline{Y_y})$
is polynomially convex in each fiber $\{y\}\times\CC^2$. \

Then, by possibly having to decrease $\e$, for all $\widehat\d>0$
there exist a family of domains $U_y\subset X_y$, $y\in B_\e$, and
a smooth map $G:\cup_{y\in B_\e}\{y\}\times\overline
U_y\rightarrow B_\e\times\CC^2$ such that the following hold for
all $y\in B_\e$:

\

(i) \ $U_y$  is homeomorphic to $Y_y$, \

(ii) \ $Y_y\subset U_y\subset Y^{\widehat\d}_y$, \

(iii) \ $d_H(U_{y_j},U_y)\rightarrow 0$  for all $y_j\rightarrow
y, \ y_j\in B_\e$, \

(iv) \ $G(y,\cdot)$  is a holomorphic embedding of $U_y$  into
$\{y\}\times\CC^2$, \

(v) \ $G(y,\overline U_y)$  satisfies the conditions in Theorem 5
when regarded as an embedded Riemann surface in the fiber
$\{y\}\times\CC^2$.

\end{proposition}
\begin{proof}
We will prove the result in the case that each fiber $Y_y$  is a
closed Riemann surface with a single component removed. We will
make some comments along the way as regards the general case,
which is essentially the same. \

We may assume that $F(\overline{Y^\d})\subset
B_\e\times\triangle\times\CC$. For any $0<r<\widehat\d$ let
$s_r\subset\overline\triangle$ denote the curve
$s_r:=\{z\in\CC;\mathrm{Im}(z)=0,-1\leq\mathrm{Re}(z)\leq -1+r\}$,
and let $S_r\subset B_\e\times\overline\triangle$  denote the
manifold $S_r:=\cup_{y\in B_\e} \{y\}\times s_r$.  Then
$f_1(S_r)\subset X$  is a smooth manifold attached to the boundary
of $Y$  with $f_1(S_r)\subset Y^{\widehat\d}\setminus Y$.  In each
fiber $Y^\d_y$  we have that $c_y:=f_1(S_r)\cap Y^\d_y$  is a
smooth curve attached to the Riemann surface $Y_y$.  \

Let $H$  denote the composition $F\circ f_1$, and let $E_r$
 denote $H(S_r)$.  Then $E_r$  is a submanifold of $B_\e\times\CC^2$, and
each fiber slice $\gamma_y:=E_r\cap(\{y\}\times\CC^2)$  is a
smooth curve attached to the embedded Riemann surface $F(Y_y)$.  \

Let us first concentrate on some fiber over $y\in B_\e$  and
explain how we can modify $F|_{Y_y^\d}$  to get all claims in the
theorem, except of course $(iii)$, for that particular fiber. The
idea is the following: We find a neighborhood $W_y$  of
$F(Y_y)\cup\gamma_y$ in $\{y\}\times\CC^2$  and an injective
holomorphic map $\psi_y:W_y\rightarrow\{y\}\times\CC^2$  such that
$\psi$  is close to the identity on $F(Y_y)$  and such that
$\psi_y$ stretches the curve $\gamma_y$  so that
$\psi_y(\gamma_y)$ intersects the cylinder
$\{y\}\times\partial\triangle\times\CC^2$ transversally and at a
single point.  For a small $\mu>0$  let $V_y^\mu$  denote the
$\mu$-neighborhood
$$
(*) \ V_y^\mu:=\{x\in Y_y^\d;d(x,Y_y\cup c_y)<\mu\}
$$
of $Y_y\cup c_y$  in $Y_y^\d$. We find a pair $(G_y,U_y)$ as in
the proposition by defining $G_y:=\psi_y\circ F$ and then
$$
(**) \
U_y:=G_y^{-1}(G_y(V_y^\mu)\cap(\{y\}\times\triangle\times\CC)).
$$
(Meaning that $U_y$  is the connected component of the pullback
that contains $Y_y$).  In the general case we attach disjoint
curves in a similar manner, one for each boundary component, and
stretch each curve.  \

More detailed we carry out the construction (still focusing on a
particular fiber) as follows: Let $m_y$ be a smoothly embedded
curve $m_y:[0,1]\rightarrow\{y\}\times\CC^2$ such that

\

(i) \ $m_y\cap F(Y^{\d}_y\setminus Y_y)\supset\gamma_y$, \

(ii) \ $(m_y\setminus\gamma_y)\cap F(\overline Y_y)=\emptyset$, \

(iii) \ The intersection
$\gamma_y\cap(\{y\}\times\partial\triangle\times\CC)$  consists of
a single point (which is not the end point), and the intersection
is transversal.

\

Let $x_0\in (0,1)$  and let $g:[0,\infty)\times [0,1]\rightarrow
[0,1]$ be an isotopy of diffeomorphisms such that

\

(a) \ $g(t,x)=x$  for all $x\in [0,x_0], t\in [0,\infty)$, \

(b) \ $\mathrm{lim}_{t\rightarrow\infty} g(t,x)=1$  for all
$x>x_0$.

\

Define an isotopy $\phi_y:[0,1]\times m_y\rightarrow m_y$  by
$\phi_y(t,x):=m_y\circ g(t,m_y^{-1}(x))$.  If $N_y$  is a small
neighborhood of $F(Y_y)$  in $\{y\}\times\CC^2$  we may define an
isotopy of diffeomorphisms $\xi_y:[0,1]\times
N_y\cup\gamma_y\rightarrow N_y\cup\gamma_y$  by
$$
\xi_y|_{N_y}:=\mathrm{Id}, \ \xi_y(t,x):=\phi(t,x) \ \mathrm{for}
\ x\in m_y.
$$
We will argue in a moment that for arbitrarily small $x_0$  and
arbitrarily large $t_0$  there is a neighborhood $W_y$  of
$F(Y_y)\cup m_y$  in $\{y\}\times\CC^2$  such we can approximate
the map $\xi_y(t_0,\cdot)$  good in $\mathcal{C}^1$-norm on
$F(Y_y)\cup m_y$  by an injective holomorphic map
$$
\psi_y:W_y\rightarrow\{y\}\times\CC^2.
$$
Granted the existence of this approximation this proves, by the
construction $(*)$  and $(**)$ above, the result (except $(iii)$)
for any particular fiber $Y_y$. \

To get $(iii)$  we carry out this construction simultaneously for
all fibers.  By possibly having to decrease $\e$   we see that we
can find a smooth submanifold $M$  of $B_\e\times\CC^2$ such that
in each fiber we have that $m_y:=M\cap\{y\}\times\CC^2$ is a
smooth curve satisfying $(i)-(iii)$  above.  Let $D:B_\e\times
[0,1]\rightarrow M$  be a diffeomorphism.  In the general case we
attach several disjoint smooth manifolds, one for each boundary
component.  For dimension reasons this does not raise a problem. \

Let $\varphi:[0,\infty)\times B_\e\times [0,1]\rightarrow
B_\e\times [0,1]$ be the isotopy $\varphi(t,y,x)=(y,g(t,x))$, and
let $\phi:[0,\infty)\times M\rightarrow M$  be the isotopy
$\phi=D\circ\varphi\circ D^{-1}$.  \

Now regard $B_\e(\t,\x)$  as the real $\e$-ball contained in
$\CC^N$, and let $\mathcal{N}$  be a small neighborhood of $F(Y)$
in $\CC^N\times\CC^2$. Define
$\xi:[0,\infty)\times(\mathcal{N}\cup M)\rightarrow
B_\e(\t,\x)\times(\mathcal{N}\cup M)$ by
$$
\xi(t,x):=x \ \mathrm{for} \ x\in\mathcal{N}, \xi(t,x):=\phi(t,x)
\ \mathrm{for} \ x\in M.
$$
Since each $F(Y_y)$  is polynomially convex in the fiber over
$\{y\}$  it follows by \cite{Sb}  that each $F(Y_y)\cup m_y$  is
polynomially convex in the fiber.  And so since
$B_\e(\t,\x)\subset\CC^N$ is real it follows that $F(Y)\cup M$  is
polynomially convex in $\CC^N$. \

By \cite{fl} we have then that for any fixed $t_0$  and $x_0$
there is a neighborhood $W$ of $F(Y)\cup M$ such that
$\xi(t_0,\cdot)$ can by approximated arbitrarily good by an
injective holomorphic map $\psi:W \rightarrow B_\e\times\CC^2$
preserving fibers, and the approximation is good in
$\mathcal{C}^1$-norm. \

Define $G:=\psi\circ F$, chose a small $\mu>0$  and define domains
$U_y$  as in $(*)$  and $(**)$  above.  If $\mu$  is small enough
then $(iii)$  follows by transversality.

\end{proof}

To prove Proposition \ref{mainmap} then, we have to construct
manifolds $X,Y$ and $Y^\d$  as above with subsets of tori as
fibers, construct a suitable map $F$, and then apply Proposition
\ref{conper}. \

Recall the Weierstrass p-function (depending on $\l$):
$$
\varrho_\l(z)=\frac{1}{z^2}+\sum_{(m,n)\in\ZZ^2\setminus
(0,0)}\frac{1}{(z-(m+n\cdot\l))^2}- \frac{1}{(m+n\cdot\l)^2}.
$$
This a meromorphic function in $z$  respecting the relation
$\sim_\l$.  Fix a 1-domain $\O(\t,0,r)$, an $\e>0$, and let
$W_\e:=\cup_{\l\in\triangle_\e(\t)}\{\l\}\times\O(\l,0,r)$. If
$\e>0$ is small enough and $p$ is close to the origin we may
define a map
$$
\widehat\phi_p(\l,z)=(\varrho_\l(z-p),\varrho_\l(z)),
$$
from $W_\e$  into $\CC^2$.

\begin{lemma}\label{jac}
For sufficiently small $\e$  and $p$  we have that
$\widehat\phi_p$  is holomorphic in the variables $(\l,z)$.   For
each fixed $\l$  we have that $\widehat\phi_p(\l,\cdot)$  embeds
$\R(\l,\O(\l,0,r))$  into $\CC^2$.
\end{lemma}
\begin{proof}
If $\e$  and $p$  is chosen small enough we have that
$\widehat\phi_{p}(\l,z)$  is holomorphic in the $z$-variable for
all fixed $\l\in\triangle_{\e}(\t)$.  To prove that $\widehat\phi$
is holomorphic in both variables we inspect the standard proof of
the fact that $\varrho_{\l}(z)$  converges as a function in the
$z$-variable. \

Following Ahlfors \cite{ah} we have for $2|z|\leq |m+n\t|$, that
$$
|\frac{1}{(z-(m+n\t))^2}-\frac{1}{(m+n\t)^2}\mid \leq\frac{10|
z|}{|m+n\t|^3}.
$$
So to prove that $\varrho_{\t}(z)$  converges it is enough to
prove that
$$
\sum_{(m,n)\in\ZZ^2\setminus (0,0)}\frac{1}{|m+n\t|^3}
$$
converges.  This in turn is proved by observing that there exists
a positive constant $K$  such that
$$
|m + n\t|\geq K(|m|+ |n|)
$$
for all $m,n\in\NN$, and then getting the estimate
$$
(*)  \sum_{(m,n)\in\ZZ^2\setminus (0,0)}\frac{1}{|m+n\t|^3}\leq
4K^{-3}\sum_{n=1}^{\infty}\frac{1}{n^2}<\infty.
$$
But $K$  may be chosen such that
$$
|m+n\l|\geq K(|m|+|n|)
$$
for all $\l$  close to $\t$, so the inequality $(*)$  holds as we
vary $\t$.  This shows that the sum $\varrho_{\l}(z)$  converges
uniformly on compacts in $W_{\e}$  in the variables $(\l,z)$. And
if the shift determined by $p$   is small enough we have that
$\widehat\phi_{p}$ is holomorphic on $W_{\e}$. \

In \cite{wd3}  we demonstrated that the map $z\mapsto
(\varrho_\l(z-p),\varrho_\l(z))$  is an embedding provided that
$2p$  is not contained in the lattice determined by $\l$. So all
$\phi_p(\l,\cdot)$  are fiberwise embeddings as long as $\e$  is
small, and $p$  is chosen close to the origin. \
\end{proof}

Let us now construct manifolds $X,Y$  and $Y^\d$  as above.  Fix
an $m$-domain $\O(\t,\x)$  and let $\e>0$.  We define
$X:=\cup_{(\l,\y)\in B_\e(\t,\x)}\{(\l,\y)\}\times\R(\O(\l))$ and
we let $\pi:X\rightarrow B_\e(\t,\x)$  be the obvious projection.
Let $q:B_\e(\t,\x)\times\CC\rightarrow X$  be the map defined by
the standard quotient map on each fiber -
$q(\l,\y,\z)=(\l,\y,[\z])$  where $[\z]$  denotes the equivalence
class of $\z$  in $\CC/\sim_\l$.  Then $q$ induces a
differentiable structure on $X$ such that each fiber $X_{(\l,\y)}$
is a closed Riemann surface which we equip with the complex
structure corresponding to $\l$. Let $m:X\times X\rightarrow\RR^+$
be a smooth metric that induces the topology.  \

Next let $V_\e=\cup_{(\l,\y)\in
B_\e(\t,\x)}\{(\l,\y)\}\times\O(\l,\y)$.  Then $Y:=q(V_\e)\subset
X$ is a submanifold $Y$  of $X$  as above.  This is seen by
defining $g_i:B_\e(\t,\x)\times\overline\triangle\rightarrow
B_\e(\t,\x)\times\CC$  by $g_i(\l,\y,t)=(\l,\y,z_i+t\cdot r_i)$
and $f_i=q\circ g_i$. \

To construct the map $F:Y^\d\rightarrow B_\e\times\CC^2$  we first
let $V^\d_\e$  denote the set $q^{-1}(Y^\d)$, and define a map
$$
\phi_p:V^\d_\e\rightarrow B_\e(\t,\x)\times\CC^2
$$
by $\phi_p(\l,\y,\z)=(\l,\y,\widehat\phi_p(\l,\z-z_1))$ (here
$z_1$ is a component of the fixed point $(\t,\x)$  and not a
variable). This is a well defined mapping if $\e$  and $p$ are
small enough. Now define a map
$$
\Phi:Y^\d\rightarrow B_\e(\t,\x)\times\CC^2
$$
by $\Phi(x)=\phi_p(q^{-1}(x))$  for $x\in Y^\d$.  This is well
defined because $\phi_p$  respects the relation $\sim_\l$  on
fibers, and it follows from Lemma \ref{jac}  that $\Phi$  is a
smooth mapping such that $\Phi|_{X_y}$  is an embedding for each
fiber $X_y$.  In the following proof of Proposition \ref{mainmap}
we use $\Phi$  to construct $F$: \

\emph{Proof of Proposition \ref{mainmap}:} Let $X,Y,Y^\d$  and
$\Phi$  be as just defined.  By Proposition \ref{perturb} there is
an open set $U\subset\CC^2$ and an injective holomorphic map
$\xi:U\rightarrow\CC^2$  such that
$\Phi(Y^\d_{(\t,\x)})\subset\{(\t,\x)\}\times U$, and such that
$\xi\circ\Phi(Y_{(\t,\x)})$  is polynomially convex in the fiber
$(\t,\x)\times\CC^2$.  Define
$$
\Psi: B_\e(\t,\x)\times U\rightarrow B_\e(\t,\x)\times\CC^2
$$
by $\Psi(\l,\y,w_1,w_2)=(\l,\y,\xi(w_1,w_2))$.\

If $\e$  is small enough we have that $\Psi\circ\Phi(Y_{(\l,\y)})$
is polynomially convex in the fiber $(\l,\y)\times\CC^2$  for all
$(\l,\y)$.  To see this choose a Runge and Stein domain
$N\subset\CC^2$ such that
$\Psi\circ\Phi(Y_{(\t,\x)})\subset\{(\t,\x)\}\times N$  and
$\Psi\circ\Phi(Y_{(\t,\x)}^\d)\cap\{(\t,\x)\}\times
N\subset\subset\Psi\circ\Phi(Y_{(\t,\x)}^\d)$.  If $\e$  is small
then $\Psi\circ\Phi(Y_{(\l,\y)}^\d)\cap\{(\l,\y)\}\times
N\subset\subset\Psi\circ\Phi(Y_{(\l,\y)}^\d)$  for all $(\l,\y)\in
B_\e(\t,\x)$, i.e.
$\Psi\circ\Phi(Y_{(\l,\y)}^\d)\cap\{(\l,\y)\}\times N$  is a
closed submanifold of $\{(\l,\y)\}\times N$.  So if $\e$  is small
the claim follows from Lemma 3.\

Define $F=\Psi\circ\Phi$ and the pair $(Y^\d,F)$ satisfies the
conditions in Proposition \ref{conper}. Let $G$ be as in
Proposition \ref{conper} and define $\psi(\O(\l,\y))$
 to be the $m$-domain corresponding to $U_{(\l,\y)}$.  Now
 $(i)-(v)$  guaranties that the conclusions of Proposition
 \ref{mainmap}  are satisfied.$\hfill\square$

\bibliography{biblio}

\begin{thebibliography}{10}

\bibitem{ah}
L.V. Ahlfors.
\newblock {\em Complex Analysis.}
\newblock McGraw Hill, 1966.

\bibitem{al}
H.~Alexander.
\newblock Explicit imbedding of the (punctured) disc into $\mathbb{C}^2$.
\newblock {\em Math.Helv.}, 52:439--544, 1977.

\bibitem{bs}
H.~Behnke and K.~Stein.
\newblock Entwicklung analytisher {F}unktionen auf {R}iemannschen {F}lachen.
\newblock {\em Math. Ann.}, 120:430--461, 1949.

\bibitem{eg}
Y.~Eliashberg and M.~Gromov.
\newblock Embeddings of {S}tein manifolds of dimension $n$ into the affine
  space of dimension $3n/2+1$.
\newblock {\em Ann.Math.}, 136:123--135, 1992.

\bibitem{fs2}
O.~Forster.
\newblock Plongements des vari\'{e}t\'{e}s de {S}tein.
\newblock {\em Comm.Math.Helv.}, 45:170--184, 1970.

\bibitem{fs}
O.~Forster.
\newblock {\em Lectures on {R}iemann Surfaces}.
\newblock Springer-Verlag, 1999.

\bibitem{fc2}
F.~Forstneri\v{c}.
\newblock The homotopy principle in complex analysis: A survey.
\newblock {\em Contemp. Math., Amer. Math. Soc., Providence, RI}, 332:73--99,
  2003.

\bibitem{fl}
F.~Forstneri\v{c} and E.~L{\o}w.
\newblock Global holomorphic equivalence of smooth manifolds in $\mathbb{C}^k$.
\newblock {\em Indiana Univ.Math.J.}, 46:133--153, 1997.

\bibitem{gs}
J.~Globevnik and B.~Stens{\o}nes.
\newblock Holomorphic embeddings of some planar domains into $\mathbb{C}^2$.
\newblock {\em Math. Ann.}, 303:579--597, 1995.

\bibitem{gz}
G.M. Goluzin.
\newblock {\em Geometric theorey of functions of a complex variable.}
\newblock American mathematical society, Providence, R.I., 1969.

\bibitem{gr}
R.C. Gunning and Rossi H.
\newblock {\em Analytic functions of several complex variables}.
\newblock Prentice-Hall, Inc., 1965.

\bibitem{hs}
Z-X. He and O.~Schramm.
\newblock Fixed points, {K}oebe uniformization, and circle packings.
\newblock {\em Ann.Math.}, 137:369--406, 1993.

\bibitem{kn}
K.~Kasahara and T.~Nishino.
\newblock As announced in math reviews.
\newblock {\em Math.Reviews.}, 38, 1969.

\bibitem{la}
H.B. Laufer.
\newblock Imbedding annuli in $\mathbb{C}^2$.
\newblock {\em J.d'Analyse Math.}, 26:187--215, 1973.

\bibitem{ma}
B.~Malgrange.
\newblock Existence et approximation des solutions des \'{e}quations aux
  d\'{e}riv\'{e}es partielles et des \'{e}quations de convolution.
\newblock {\em Ann. Inst. Fourier}, 6:271--354, 1955-56.

\bibitem{sc}
J.~Schurmann.
\newblock Embeddings of {S}tein spaces into affine spaces of minimal dimension.
\newblock {\em Math.Ann.}, 307:381--399, 1997.

\bibitem{Sb}
G.~Stolzenberg.
\newblock Uniform approximation on smooth curves.
\newblock {\em Acta Math.}, 115:185--198, 1966.

\bibitem{cf}
M.~\v{C}erne and F.~Forstneri\v{c}.
\newblock Embedding some bordered {R}iemann surfaces in the affine plane.
\newblock {\em Math. Res. Lett.}, 9:683--696, 2002.

\bibitem{wd3}
E.~F. Wold.
\newblock Embedding {R}iemann surfaces into $\mathbb{C}^2$.
\newblock {\em {I}nternat.{J}.{M}ath}, 17:963--974, 2006.

\bibitem{wd2}
E.~F. Wold.
\newblock Proper holomorphic embeddings of finitely and some infinitely
  connected subsets of $\mathbb{C}$ into $\mathbb{C}^2$.
\newblock {\em Math.Z.}, 252:1--9, 2006.

\end{thebibliography}

\end{document}